\documentclass[10pt]{article}
\usepackage{latexsym}
\usepackage{indentfirst}
\usepackage{amssymb,amsfonts,amsmath,amsthm}
\usepackage{graphicx}
\usepackage{mathrsfs}
\usepackage{array}
\usepackage{color}
\usepackage{epstopdf}
\usepackage[all]{xy}
\usepackage{hyperref}
\usepackage{url}

\usepackage{geometry}
\newtheorem{thm}{Theorem}[section]
\newtheorem{prop}[thm]{Proposition}
\newtheorem{lem}[thm]{Lemma}
\newtheorem{conj}[thm]{Conjecture}
\newtheorem{exam}[thm]{Example}

\newtheorem{rmk}[thm]{Remark}
\newtheorem{dfn}[thm]{Definition}

\numberwithin{equation}{section}

\newcommand{\frakm}{{\mathfrak m}}

\newcommand{\frakG}{{\mathfrak G}}

\newcommand{\frakM}{{\mathfrak M}}
\newcommand{\frakN}{{\mathfrak N}}

\newcommand{\frakS}{{\mathfrak S}}

\newcommand{\bE}{{\mathbb E}}
\newcommand{\bF}{{\mathbb F}}

\newcommand{\bM}{{\mathbb M}}
\newcommand{\bN}{{\mathbb N}}

\newcommand{\bQ}{{\mathbb Q}}

\newcommand{\bV}{{\mathbb V}}

\newcommand{\bZ}{{\mathbb Z}}

\newcommand{\calE}{{\mathcal E}}
\newcommand{\calF}{{\mathcal F}}
\newcommand{\calG}{{\mathcal G}}
\newcommand{\calH}{{\mathcal H}}
\newcommand{\calI}{{\mathcal I}}

\newcommand{\calL}{{\mathcal L}}
\newcommand{\calM}{{\mathcal M}}

\newcommand{\calO}{{\mathcal O}}

\newcommand{\calV}{{\mathcal V}}

\newcommand{\calX}{{\mathcal X}}

\newcommand{\Coker}{{\mathrm{Coker}}}       
\newcommand{\End}{{\mathrm{End}}}           
\newcommand{\Ext}{{\mathrm{Ext}}}           
\newcommand{\Gal}{{\mathrm{Gal}}}           
\newcommand{\Hom}{{\mathrm{Hom}}}           
\newcommand{\Ima}{{\mathrm{Im}}}            
\newcommand{\Ker}{{\mathrm{Ker}}}           

\newcommand{\Mod}{{\mathrm{Mod}}}           

\newcommand{\Rep}{{\mathrm{Rep}}}           

\newcommand{\Spf}{{\mathrm{Spf}}}           
\newcommand{\Spec}{{\mathrm{Spec}}}         
\newcommand{\Vect}{{\mathrm{Vect}}}         

\newcommand{\GL}{{\mathrm{GL}}}             

\newcommand{\cris}{{\mathrm{cris}}}         
\newcommand{\et}{{\mathrm{\text{\'e}t}}}    

\newcommand{\tor}{{\mathrm{tor}}}           
 
\newcommand{\Art}{{\mathrm{Art}}}  

\usepackage{relsize}
\usepackage[bbgreekl]{mathbbol}
\DeclareSymbolFontAlphabet{\mathbb}{AMSb} 
\DeclareSymbolFontAlphabet{\mathbbl}{bbold}
\newcommand{\Prism}{{\mathlarger{\mathbbl{\Delta}}}} 

\begin{document}
\title{Prismatic construction of crystalline deformation rings up to $p$-power torsion}
\date{}
\author{Yu Min}
\maketitle

\begin{abstract}
    We study the deformation of Galois representations which are \'etale realizations of finitely generated prismatic $F$-crystals and show the deformation functor is prorepresented by a complete Noetherian local ring whose $\bZ_p$-flat quotient coincides with Kisin's crystalline deformation ring. As a consequence, we obtain a universal prismatic $F$-crystal with coefficients in Kisin's crystalline deformation ring.
\end{abstract}
\tableofcontents

\section{Introduction}

The crystalline deformation ring constructed by Kisin in \cite{kisin2008potentially} plays a central role in many important questions like the modularity lifting problem and the Breuil--M\'ezard conjecture. Despite its name, this ring is not obtained by a deformation theoretic construction. Instead, Kisin's construction involves taking scheme theoretic closure from the generic fiber. So it is reasonable that this deformation ring is very hard to investigate, especially its special fiber.

It is then a natural question to seek a direct deformation theoretic construction of the crystalline deformation ring. In this short note, we use the prismatic theory to propose another related deformation functor, based on \cite{BS23}. 

Let $K$ be a finite extension of $\bQ_p$ with ring of integers $\calO_K$ and  residue field $k$. Let $e$ be the absolute ramification degree of $K$ and let $G_K$ be the absolute Galois group of $K$. Fix a finite field $\bF$ of characteristic $p$ and let $\Art_{W(\bF)}$ be the category of finite Artinian local rings over $W(\bF)$ with fixed isomorphisms between their residue field and $\bF$. 

A deformation theoretic construction requires a reasonable definition of crystalline $R$-representations for any $R\in \Art_{W(\bF)}$. Using Fontaine's definition of torsion crystalline representations, Tong Liu \cite{liu2006torsion} constructed a deformation ring which differs from Kisin's by $p$-power torsion. 

Thanks to \cite{BS23}, we now have a nice description of crystalline $\bZ_p$-representations using prismatic theory. Let $(\calO_K)_{\Prism}$ be the absolute prismatic site of $\calO_K$. Consider the Breuil--Kisin prism $(\frakS=W(k)[[u]],E(u))$, where $E(u)$ is the Eisenstein polynomial of a fixed uniformiser $\pi$ of $\calO_K$. The Breuil--Kisin prism $(\frakS=W(k)[[u]],E(u))$ is a cover of the final object of the topos ${\rm Shv}((\calO_K)_{\Prism})$. Let $(\frakS^{(\bullet)},E)$ denote the \v Cech nerve corresponding to this cover. In particular, we write  $(\frakS^{(1)},E):=(\frakS,E)\times(\frakS,(E))$ and $(\frakS^{(2)},E):=(\frakS,E)\times(\frakS,(E))\times(\frakS,(E))$. Bhatt--Scholze proved the following theorem.

\begin{thm}[{\cite[Theorem 1.2]{BS23}}]
    There is an equivalence of categories
    \[
    \Vect^{\varphi}((\calO_K)_{\Prism},\calO_{\Prism})\simeq \Rep^{\cris}_{\bZ_p}(G_K)
    \]
    where the left hand side is the category of prismatic $F$-crystals of finite projective modules, and the right hand side is the category of crystalline $\bZ_p$-representations of $G_K$.
\end{thm}

Prismatic $F$-crystals of finite projective modules are equivalent to finite projective Breuil--Kisin modules equipped with descent data $(\frakM,\epsilon_{\frakM})$ relative to the \v Cech nerve $(\frakS^{(\bullet)},E)$, i.e.
\begin{enumerate}
    \item $\frakM$ is a finite projective module equipped with an isomorphism $(\varphi^*\frakM)[1/E]\cong \frakM[1/E]$.
    \item $\epsilon_{\frakM}:\frakM\otimes_{\frakS,p_0}\frakS^{(1)}\cong \frakM\otimes_{\frakS,p_1}\frakS^{(1)}$ satisfying the cocycle condition over $\frakS^{(2)}$.
\end{enumerate}

This theorem is the starting point of our deformation functor. As the crystalline deformation ring is a quotient of the unrestricted local deformation ring, its corresponding functor on $\Art_{W(\bF)}$ should be a full subfunctor of the unrestricted local deformation functor. In particular, we just need to single out the ``crystalline $R$-representations" and the isomorphisms between them should just be the isomorphisms as $R$-representations. One possible definition of crystalline $R$-representations is those $R$-representations which are the \'etale realisation of some prismatic $F$-crystals of finite projective modules over $\calO_{\Prism,R}:=\calO_{\Prism}\otimes_{\bZ_p}R$. 

Hence the first thing we need to do is to have an \'etale realisation functor with $R$-coefficients. This follows from the following theorem.

\begin{thm}[{\cite[Theorem 1.1]{Min23}}]
    There is an equivalence of groupoids
    \[
    \Vect((\calO_K)_{\Prism},(\calO_{\Prism}\otimes_{\bZ_p}R)[1/\calI_{\Prism}])^{\varphi=1}\simeq \calX_{\rm EG}(R)\simeq \Rep_R(G_K)
    \]
    for any $R\in \Art_{W(\bF)}$, where $\Vect((\calO_K)_{\Prism},(\calO_{\Prism}\otimes_{\bZ_p}R)[1/\calI_{\Prism}])^{\varphi=1}$ is the groupoid of Laurent $F$-crystals with $R$-coefficients, $\calX_{\rm EG}(R)$ is the groupoid of $R$-points of the Emerton--Gee stack, which are equivalent to the $R$-representations of $G_K$.
\end{thm}

A Laurent $F$-crystal in $\Vect((\calO_K)_{\Prism},(\calO_{\Prism}\otimes_{\bZ_p}R)[1/\calI_{\Prism}])^{\varphi=1}$ is the data $(M, \epsilon_M)$ where
\begin{enumerate}
    \item $M$ is a finite projective \'etale $\varphi$-module over $\frakS_R[1/E]:=(\frakS\otimes_{\bZ_p}R)[1/E]=(W(k)\otimes_{\bZ_p}R)((u))$.
    \item $\epsilon_M: M\otimes_{\frakS_R[1/E],p_0}\frakS^{(1)}_R[1/E]\cong M\otimes_{\frakS_R[1/E],p_1}\frakS^{(1)}_R[1/E]$ satisfying the cocycle condition over $\frakS^{(2)}_R[1/E]$.
\end{enumerate}

Regarding $R$-representations of $G_K$ as Laurent $F$-crystals with $R$-coefficients, the \'etale realisation functor just sends a prismatic $F$-crystal $(\frakM,\epsilon_{\frakM})$ to the Laurent $F$-crystal $(\frakM[1/E],\epsilon_{\frakM}[1/E])$.

But it turns out admitting a model of finite projective prismatic $F$-crystal is too restrictive to get a pro-representable deformation functor. The solution to this difficulty is inspired by the following observation.

\begin{lem}[Lemma \ref{observation}]\label{intro-observation}
     Let $\rho$ be a $\bZ_p$-representation of $G_K$ which is the \'etale realisation of a finitely generated prismatic $F$-crystal. Then $\rho$ is crystalline.
\end{lem}

Let us first explain the meaning of finitely generated prismatic $F$-crystal. 
\begin{dfn}
    By a finitely generated prismatic $F$-crystal over $\calO_{\Prism}$, we mean a pair $(\frakM,\epsilon_{\frakM})$ where
\begin{enumerate}
    \item $\frakM$ is a finitely generated torsion free Breuil--Kisin module over $\frakS$,
    \item $\epsilon_{\frakM}:\frakM\widehat\otimes_{\frakS,p_0}\frakS^{(1)}\cong \frakM\widehat\otimes_{\frakS,p_1}\frakS^{(1)}$ satisfying the cocycle condition over $\frakS^{(2)}$, where the completion on both sides is $(p,E)$-adic completion.
\end{enumerate}

We say $(\frakM,\epsilon_{\frakM})$ is of height at most $h$ if $\frakM$ is of height at most $h$.
\end{dfn}

By \cite[Lemma 3.5]{DLMS24}, we know the maps $p_i:\frakS\to \frakS^{(1)}$ (resp. $q_i:\frakS\to \frakS^{(2)}$) are all faithfully flat for $i=0,1$ (resp. $i=0,1,2$). This implies $\frakM\otimes_{\frakS,p_i}\frakS^{(1)}$ is already $(p,E)$-complete by \cite[Remark 3.12]{DLMS24}. 

\begin{dfn}
    Let $(M,\epsilon_M)$ be a Laurent $F$-crystal in $\Vect((\calO_K)_{\Prism},(\calO_{\Prism})[1/\calI_{\Prism}]^{\wedge}_p)^{\varphi=1}$. We say $(M,\epsilon_M)$ admits a finitely generated prismatic $F$-model if there is a finitely generated prismatic $F$-crystal $(\frakM,\epsilon_{\frakM})$ such that $(\frakM\otimes_{\frakS}\frakS[1/E]^{\wedge}_p,\epsilon_{\frakM}[1/E]^{\wedge}_p)$ is isomorphic to $(M,\epsilon_M)$.
\end{dfn}

We can also define the above notions with $R$-coefficients.
\begin{dfn}
   
 By a finitely generated prismatic $F$-crystal over $\calO_{\Prism,R}:=\calO_{\Prism}\otimes_{\bZ_p}R$, we mean a pair $(\frakM,\epsilon_{\frakM})$ where
\begin{enumerate}
    \item $\frakM$ is a finitely generated $u$-torsion-free Breuil--Kisin module over $\frakS_R:=\frakS\otimes_{\bZ_p}R$,
    \item $\epsilon_{\frakM}:\frakM\otimes_{\frakS_R,p_0}\frakS_R^{(1)}\cong \frakM\otimes_{\frakS_R,p_1}\frakS_R^{(1)}$ satisfying the cocycle condition over $\frakS^{(2)}_R$. Again, we do not need to take the $E$-adic completion as both sides are already $E$-complete.
\end{enumerate}
We say $(\frakM,\epsilon_{\frakM})$ is of height at most $h$ if $\frakM$ is of height at most $h$.
\end{dfn}

\begin{dfn}
    Let $(M,\epsilon_M)$ be a Laurent $F$-crystal in $\Vect((\calO_K)_{\Prism},(\calO_{\Prism,R})[1/\calI_{\Prism}])^{\varphi=1}$. We say $(M,\epsilon_M)$ admits a finitely generated prismatic model if there is a finitely generated prismatic $F$-crystal $(\frakM,\epsilon_{\frakM})$  such that $(\frakM\otimes_{\frakS_R}\frakS_R[1/E],\epsilon_{\frakM}[1/E])$ is isomorphic to $(M,\epsilon_M)$.
 
\end{dfn}

Now let $\bar\rho:G_K\to \GL_n(\bF)$ be an $\bF$-representation. We assume $\bar\rho$ admits a finitely generated prismatic $F$-model. By virtue of Lemma \ref{intro-observation}, we will consider those deformations which admit finitely generated prismatic $F$-models of height at most $h$. It turns out this deformation ring is exactly the crystalline deformation ring after quotientIng out all the $p$-power torsion.

\begin{dfn}
   Let $\calF^{h}:\Art_{W(\bF)}\to Groupoids$ be the functor sending $R$ to the groupoid of $R$-representations which admit a finitely generated prismatic $F$-model of height $h$. Define the functor $\calF^{h}_{\bar\rho}:\Art_{W(\bF)}\to Groupoids$ as follows
   \[
   \calF^{h}_{\bar\rho}(R):=\calF^{h}(R)\times_{\calF^{h}(\bF)}*_{\bar\rho}.
   \]
We can also define a framed version $\calF^{h,\Box}_{\bar\rho}$ by sending $R\in \Art_{W(\bF)}$ to the set $\{\rho\in \GL_d(R):\rho\equiv \bar\rho\mod (\frakm_R)\ \text{and $\rho$ admits a finitely generated prismatic $F$-model of height $h$}\}$.
\end{dfn}
Let $R_{\bar\rho}^{\Box}$ denote the framed unrestricted local deformation ring and $R_{\bar \rho}^{\cris,h,\Box}$ be Kisin's crystalline deformation ring for a fixed height $h$. Our main theorem is as follows.

\begin{thm}[{Theorem \ref{main1}, \ref{main2}, \ref{universal prismatic crystal}}]
\begin{enumerate}
    \item 
    The functor $\calF^{h,\Box}_{\bar\rho}$ is pro-represented by a complete Noetherian local ring $S^{h,\Box}_{\bar \rho}$, which is a quotient of $R_{\bar\rho}^{\Box}$. 
    \item The $\bZ_p$-flat quotient of $S^{h,\Box}_{\bar \rho}$ is exactly Kisin's crystalline deformation ring, i.e. $S^{h,\Box}_{\bar \rho}/(S^{h,\Box}_{\bar \rho}[p^\infty])=R_{\bar \rho}^{\cris,h,\Box}$.
    \item The universal crystalline Galois deformation $G_K\to \GL_d(R_{\bar \rho}^{\cris,h,\Box})$ admits a maximal finitely generated prismatic $F$-model of height $h$.
    \end{enumerate}
\end{thm}

\begin{rmk}
All results can be easily generalised to potentially crystalline deformation rings. Indeed, to study representations which become crystalline over a fixed extension $L/K$, we just need to consider prismatic $F$-crystals over $(\calO_L)_{\Prism}$ in our constructed deformation functor. By using log-prismatic theory, we can also obtain a semistable version.
\end{rmk}

\section*{Akknowledgments}
We would like to thank Toby Gee, Eugen Hellmann, Yongquan Hu, Vytautas Paškūnas for helpful discussions. This work is partially supported by the European Research Council (ERC) under the European Union’s
Horizon 2020 research and innovation programme (grant agreement No. 884596) and the Hong Kong RGC ECS
grant 26310126.

\section*{AI disclosure}
The early version of this work was obtained without use of AI.  Later, ChatGPT pointed out a gap in our previous proof of Theorem \ref{universal prismatic crystal} and helped fix the gap. It also suggested the argument proving Lemma \ref{vanishing} in the current version. We are responsible for the correctness of all these arguements.
\section{Deformation functor}

Now we check that $\calF^{h}_{\bar\rho}$ satisfies the Rim-Schlessinger conditions. Recall that a surjective map $B_1\twoheadrightarrow B_2$ in $\Art_{\bZ_p}$ is a small extension if its kernel is principally generated and is killed by $\frakm_{B_1}$.

\begin{lem}\label{Rim-Schlessinger}
    For any pullback diagram in $Art_{W(\bF)}$
\begin{equation*}
    \xymatrix{
    \tilde A:=A_1\times_AA_2\ar[r]\ar[d]& A_2\ar[d]\\
    A_1\ar[r]& A
    }
\end{equation*}
such that $A_2\to A$ is a small extension, the induced functor
\[
\iota: \calF^{h}_{\bar\rho}(\tilde A)\to \calF^{h}_{\bar\rho}(A_1)\times_{\calF^{h}_{\bar\rho}(A)}\calF^{h}_{\bar\rho}(A_2)
\]
is an equivalence.
\end{lem}

\begin{proof}
By the definition of the isomorphisms in $\calF^{h}$, we see that $\iota$ is fully faithful as this is true for the deformation of Galois representations.

Now we focus on proving that $\iota$ is essentially surjective. Let $(\bV_1,\bV,\bV_2)$ be an object in the fiber product of
\[
\calF^{h}_{\bar\rho}(A_1)\to \calF^{h}_{\bar\rho}(A)\leftarrow \calF^{h}_{\bar\rho}(A_2).
\]
Then $\bV\cong \bV_1\otimes_{A_1}A$ and $\bV\cong \bV_2\otimes_{A_2}A$. Note that there is a pullback diagram of representations
\begin{equation*}
    \xymatrix{
    \calV\ar[r]\ar[d]&\bV_2\ar[d]\\
    \bV_1\ar[r]&\bV
    }
\end{equation*}
as the local deformation functor is representable by $R_{\bar\rho}$. The key point is that we need to find a finitely generated prismatic model of height at most $h$ of $\calV$.

Let $\bM_1,\bM_2$ be finitely generated prismatic models of height at most $h$ of $\bV_1$ and $\bV_2$ respectively. Then the image of $\bM_2\otimes_{A_2}A$ in $\bM_2[1/u]\otimes_{A_2}A$ is a finitely generated prismatic model of $\bV$. To see this, note that $\bM_2$ corresponds to a finitely generated $u$-torsion free Breuil--Kisin module $\frakM_2$ over $\frakS_{A_2}$ together with an $\frakS_{A_2}^{(1)}$-linear isomorphism $\epsilon: \frakM_2\otimes_{\frakS_{A_2},p_1}\frakS_{A_2}^{(1)}\cong \frakM_2\otimes_{\frakS_{A_2},p_2}\frakS_{A_2}^{(1)}$ satisfying the cocycle condition. By modulo $I=\Ker(A_2\twoheadrightarrow A)$, we get a finitely generated Breuil--Kisin module $\frakM_2/I$ over $\frakS_A$ together with the induced descent data, but not necessarily $u$-torsion free. Define $\frakM:=(\frakM_2/I)/(\frakM_2/I[u^{\infty}])$. Then $\frakM$ is a finitely generated $u$-torsion free Breuil--Kisin module equipped with the induced descent data, as both $p_0,p_1: \frakS_{A}\to \frakS_A^{(1)}$ are flat. Let $\bM$ be its corresponding prismatic $F$-crystal, which then is a finitely generated model of $\bV$. It is not clear whether $\bM$ is of height at most $h$.

Note that there exists an integer $n$ such that the map $\frakM_1\to \frakM[1/u]$ factors through $\frac{\frakM}{u^n}$. Now we consider the fiber product
\begin{equation*}
    \xymatrix{
    \tilde \frakM\ar[r]\ar[d]& \frac{\frakM_2}{u^n}\ar[d]\\
    \frakM_1\ar[r]& \frac{\frakM}{u^n}.
    }
\end{equation*}
We claim that $\tilde\frakM$ is a finitely generated $\frakS_{\tilde A}$-module. Note that $\frakM_1$ is a finitely generated module over $\frakS_{A_1}$. Let $e_1,\cdots, e_t$ be a set of generators of $\frakM_1$. Let $x_1=a_1e_1+\cdots +a_te_t$ be an arbitrary element in $\frakM_1$. As both $A_2\to A$ and $\frakM_2\to \frakM$ are surjective, there exists an element $(x_1,y_1)=\sum(a_i,b_i)(e_i,\tilde e_i)$ where $b_i$ is any lift in $\frakS_{A_2}$ of the image of $a_i$ in $\frakS_A$ and $\tilde e_i$ is any lift in $\frac{\frakM_2}{u^n}$ of the image of $e_i$ in $\frac{\frakM}{u^n}$. Let $N$ be the submodule of $\tilde\frakM$ generated by the elements $(e_i,\tilde e_i)$. Then any element $(x,y)\in \tilde \frakM$ can be written as $(x,y)=(x,\tilde y)+(0,z)$ such that $(x,\tilde y)\in N$ and $z\in \ker(\frac{\frakM_2}{u^n}\to \frac{\frakM}{u^n})$.

Note that \[\tilde N:=\ker(\frac{\frakM_2}{u^n}\to \frac{\frakM}{u^n})=\frac{\Ker(\frakM_2\twoheadrightarrow\frakM)}{u^n}.\]
As $\tilde N$ is a submodule of $\frac{\frakM_2}{u^n}$, it must be finitely generated over $\frakS_{A_2}$. Note that $\frakm_{A_2}I=0$ by the definition of the small extension. Then we see $\frakm_{A_2}\tilde N=0$. So $\tilde N$ is indeed a finitely generated $k[[u]]$-module. Let $\alpha_1,\cdots,\alpha_s$ be a set of generators of $\tilde N$. Then any element $z\in \tilde N$ can be written as $z=z_1\alpha_1+\cdots+ z_s\alpha_s$ with $z_i\in \frakS_{A_2}$. Now as the residue fields of $A_1,A,A_2$ are all $\bF_p$, we have
\[
\frakS_{A_1}/\frakm_{A_1}=\frakS_{A}/\frakm_{A}=\frakS_{A_2}/\frakm_{A_2}=k[[u]].
\]

Let $\bar z_i$ be the image of $z_i$ in $\frakS_{A_2}/\frakm_{A_2}=k[[u]]$ and $\beta_i$ be any lift in $\frakS_{A_1}$ of $\bar z_i$. Let $\tilde z_i$ be any lift in $\frakS_{A_2}$ of the image in $\frakS_A$ of $\beta_i$. As $\frakm_{A_2}\tilde N=0$, we see
\[
z=\sum z_i\alpha_i=\sum\tilde z_i\alpha_i.
\]
So $(0,z)=\sum (\beta_i,\tilde z_i)(0,\alpha_i)$, which means $(0,\tilde N)$ is also finitely generated over $\frakS_{\tilde A}$. So $\tilde M=N+(0,\tilde N)$ is also finitely generated over $\frakS_{\tilde A}$. Since both $\frakM_1$ and $\frac{\frakM_2}{u^n}$ are $u$-torsion free, so is $\tilde \frakM$. Now consider the pullback diagram
\begin{equation*}
    \xymatrix{
    \calM\ar[r]\ar[d]& \frakM_2\ar[d]\\
    \tilde \frakM\ar[r]& \frac{\frakM_2}{u^n}.
    }
\end{equation*}
As $\frakM_2\to \frac{\frakM_2}{u^n}$ is injective, so is $\calM\to \tilde \frakM$. So $\calM$ is also a finitely generated $u$-torsion-free $\frakS_{\tilde A}$-module. Note that $\calM$ also fits into the following pullback diagram
\begin{equation*}
    \xymatrix{
    \calM\ar[r]\ar[d]& \frakM_2\ar[d]\\
     \frakM_1\ar[r]& \frakM[1/u].
    }
\end{equation*}
So $\calM$ does not depend on $n$. Since $p_0,p_1:\frakS\to \frakS^{(1)}$ are both flat, we can actually obtain that $\calM$ is equipped with descent data $\epsilon_{\calM}:\calM\otimes_{\frakS_{\tilde A},p_0}\frakS^{(1)}_{\tilde A}\cong \calM\otimes_{\frakS_{\tilde A},p_1}\frakS^{(1)}_{\tilde A}$ by Lazard's theorem, i.e. any flat module can be written as a filtered colimit of finite free modules.

Now note that the Frobenius $\varphi:\frakS\to \frakS$ is finite flat. We also have the pullback diagram
\begin{equation*}
    \xymatrix{
    \varphi^*\calM\ar[r]\ar[d]& \varphi^*\frakM_2\ar[d]\\
    \varphi^*\tilde \frakM\ar[r]& \varphi^*\frac{\frakM}{u^n}.
    }
\end{equation*}
 So $\calM$ is a finitely generated $u$-torsion-free Breuil--Kisin module. This implies $\calM$ is a finitely generated prismatic model of $\calV$.

Note that there is an injective map between the pullback diagrams
\begin{equation*}
    \xymatrix{
    E^h\calM\ar[r]\ar[d]& E^h\frakM_2\ar[d]\\
    E^h\frakM_1\ar[r]&\frac{\frakM}{u^n}[1/u]
    }
\end{equation*}
and 
\begin{equation*}
    \xymatrix{
    \varphi^*\calM\ar[r]\ar[d]& \varphi^*\frakM_2\ar[d]\\
    \varphi^*\frakM_1\ar[r]&\varphi^*\frac{\frakM}{u^n}[1/u].
    }
\end{equation*}
The latter also admits an injective map into
\begin{equation*}
    \xymatrix{
    \calM\ar[r]\ar[d]& \frakM_2\ar[d]\\
    \frakM_1\ar[r]&\frac{\frakM}{u^n}[1/u].
    }
\end{equation*}
    So we get $E^h\calM\subset \varphi^*\calM\subset \calM$, i.e. $\calM$ is a finitely generated prismatic model of height at most $h$ of $\calV$. The functor $\iota$ is then an equivalence.
\end{proof}

Recall the following theorem of Schlessinger\cite[\href{https://stacks.math.columbia.edu/tag/06JM}{Tag 06JM}]{stacks-project} on prorepresentable functors.

\begin{thm}
    Let $\calH:\Art_{W(\bF)}\to Sets$ be a functor. Then $\calH$ is prorepresentable if and only if $(a)$ $\calH$ is a deformation functor, i.e. $\calH(\bF)=\{*\}$ and $\calH$ satisfies the Rim-Schlessinger condition, $(b)$ $\dim_{\bF}T\calH$ is finite, and $(c)$ the map $Der_{W(\bF)}(\bF,\bF)\to T\calH$ is injective (which is trivial in this case).
\end{thm}

Using the above theorem, we can show $\calF_{\bar\rho}^{h,\Box}$ is prorepresentable.
\begin{thm}\label{main1}
    The functor $\calF_{\bar\rho}^{h,\Box}$ is prorepresentable by a complete Noetherian local ring $S^{h,\Box}_{\bar \rho}$, which is a quotient of $R_{\bar\rho}^{\Box}$. 
\end{thm}
\begin{proof}
    Lemma \ref{Rim-Schlessinger} says $\calF_{\bar\rho}^{h,\Box}$ satisfies the Rim-Schlessinger condition. Note that $\calF_{\bar\rho}^{h,\Box}$ is a full subfunctor of the framed local deformation functor $\Spf(R_{\bar \rho}^{\Box})$. The latter satisfies conditions $(b)$ and $(c)$, which then implies $\calF_{\bar\rho}^{h,\Box}$ also satisfies these conditions. So we are done.
\end{proof}

 Let $\calO_B$ be any finite flat $W(\bF)$-algebra. For any ring map $R_{\bar\rho}^{\Box}\to \calO_B$, we get an induced Galois representation $\rho_B:G_K\to \GL_d(R_{\bar\rho}^{\Box})\to\GL_d(\calO_B)$. Then we have the following proposition.
 \begin{prop}\label{crys-factorisation}
     If $\rho_B$ is a crystalline representation with Hodge--Tate weights in $[0,h]$, then the ring map $R_{\bar\rho}^{\Box}\to \calO_B$ factors through $R_{\bar\rho}^{\cris,\Box}\to \calO_B$.
 \end{prop}
 \begin{proof}
     Let $A$ be the image of the map $R_{\bar\rho}^{\Box}\to \calO_B$. Then $A$ is a $\bZ_p$-flat complete Noetherian local ring with residue field $\bF$. In particular, the induced representation $\rho_A: G_K\to \GL_d(R_{\bar\rho}^{\Box})\to\GL_d(A)$ is a deformation of $\bar\rho$. Note that $\rho_A[1/p]$ is a $\bQ_p$-subrepresentation of $\rho_B[1/p]$, hence is a crystalline $\bQ_p$-representation with Hodge--Tate weights in $[0,h]$. 

     Regard $\rho_B$ as a $\bQ_p$-vector space of dimension $nd$. where $n=\dim_{\bQ_p}\calO_B[1/p]$. Then there is a rank $nd$ prismatic $F$-crystal $\bM$ in $\Vect^{\varphi}((\calO_K)_{\Prism},\calO_{\Prism})$ corresponding to $\rho_B$. Concretely, $\bE$ consists of the following data
    \begin{enumerate}
        \item $\frakM=\bM(\frakS,E)$ is a Breuil--Kisin module of dimension $nd$,
        \item $\epsilon:\frakM\otimes_{\frakS,p_0}\frakS^1\cong \frakM\otimes_{\frakS,p_1}\frakS^1$ which satisfies the cocycle condition.
    \end{enumerate}

    Note we have an injection $A\hookrightarrow \End_{\bZ_p[G_K]}(\rho_B)$. As the functor
    \[
    \Rep_{\bZ_p}^{\cris}(G_K)\to \Vect^{\varphi}(\frakS)
    \]
    is also fully faithful by \cite[Theorem 0.1]{kisin2006crystalline} or \cite[Theorem 1.3]{BS23}, we see that $\frakM$ is indeed a finitely generated $u$-torsion free Breuil--Kisin module of height $h$ and $\bM$ is a finitely generated prismatic $F$-model of height $h$ of $\rho_B$. So we get the factorization $R_{\bar\rho}^{\Box}\to S_{\bar\rho}^{\Box}\to A\to \calO_B$. 
 \end{proof}

We claim the converse is also true.

\begin{prop}\label{crystalline}
  Let $\calO_B$ be any finite flat $W(\bF)$-algebra. For any map $S_{\bar\rho}^{\Box}\to \calO_B$, the induced Galois representations $G_K\to \GL_d(\calO_B)$ is a crystalline representation with Hodge--Tate weights in $[0,h]$.
\end{prop}

Before proceeding to the proof, we first present an important observation.

\begin{prop}[{\cite[Proposition 4.3]{BMS18}}]\label{key prop}
    Let $(\frakM,\varphi_{\frakM})$ be a finitely generated torsion-free Breuil--Kisin module of height $h$. Then there is a canonical exact sequence of Breuil--Kisin modules
    \[
    0\to (\frakM,\varphi_{\frakM})\to (\frakM_{\rm free},\varphi_{\frakM_{\rm free}})\to (\bar\frakM,\varphi_{\bar\frakM})\to 0.
    \]
    In particular, $\bar\frakM$ is killed by some power of $(p,u)$ and $\frakM_{\rm free}$ is of height $h$.
\end{prop}
\begin{proof}
We only need to prove $\frakM_{\rm free}$ is of height $h$ as the rest is already proved in \cite[Proposition 4.3]{BMS18}. Recall that $\frakM_{\rm free}$ is given by $H^0(\Spec(\frakS)\setminus\{\frakm\},\calE)$ where $\calE$ is the sheaf corresponding to $\frakM$. This construction is functorial. And $E^h\frakM\subset \varphi^*\frakM\subset \frakM$ correspond to $E^h\frakM_{\rm free}\subset\varphi^*\frakM_{\rm free}\subset \frakM_{\rm free}$ under this construction, which shows $\frakM_{\rm free}$ is of height $h$.

\end{proof}

\begin{lem}\label{observation}
   Let $\rho$ be a $\bZ_p$-representation admitting a finitely generated prismatic model of height $h$. Then $\rho$ is crystalline with Hodge--Tate weights in $[0,h]$.
\end{lem}
\begin{proof}
    Let $(\frakM,\epsilon:\frakM\otimes_{\frakS,p_0}\frakS^{(1)}\cong \frakM\otimes_{\frakS,p_1}\frakS^{(1)})$ be a finitely generated prismatic model of $\rho$. Then $\frakM$ is a torsion-free finitely generated Breuil--Kisin module.

    By Proposition \ref{key prop}, we have a finite free Breuil--Kisin module $\frakM_{\rm free}$ of height $h$ corresponding to $\frakM$. Note that we have
    \[\frakM_{\rm free}=\frakM_{\rm free}[1/p]\cap \frakM_{\rm free}[1/E]^{\wedge}_p=\frakM[1/p]\cap \frakM[1/E]^{\wedge}_p\subset \frakM[1/E]^{\wedge}_p[1/p].
    \]

    As $p_0,p_1:\frakS\to \frakS^{(1)}$ are both faithfully flat, we also have
    \[
    \frakM_{\rm free}\otimes_{\frakS,p_i}\frakS^{(1)}=\frakM[1/p]\otimes_{\frakS,p_i}\frakS^{(1)}\cap \frakM[1/E]^{\wedge}_p\otimes_{\frakS,p_i}\frakS^{(1)}.
    \]
    
    So we get an isomorphism $\epsilon_{\rm free}: \frakM_{\rm free}\otimes_{\frakS,p_0}\frakS^{(1)}\cong \frakM_{\rm free}\otimes_{\frakS,p_1}\frakS^{(1)}$. This means we get a finite free prismatic $F$-crystal $(\frakM_{\rm free},\epsilon_{\rm free})$. In particular, $(\frakM,\epsilon)$ and $(\frakM_{\rm free},\epsilon_{\rm free})$ induce the same Laurent $F$-crystal, i.e. the same $\bZ_p$-representation. So $\rho$ is crystalline.
\end{proof}

The above lemma saves us from sticking to being finite projective. Now we record the following finiteness result due to Hellmann.

\begin{lem}[{\cite[Lemma 4.1]{Hel12}}]\label{torsion bounded}
    Let $(M,\varphi_M)$ be a finite free \'etale $\varphi$-module over $W(k)/p^{n+1}((u))$ of rank $d$. Let $(\frakM,\varphi_{\frakM})$ be a finitely generated Breuil--Kisin module of height $h$, which generates $(M,\varphi_M)$. Then the $k$-dimension of the $u^{\infty}$-torsion part of the finitely generated $k[[u]]$-module $\frakM/p\frakM$ is bounded by
    \[
    \dim_{k}(\frakM/p\frakM)_{\tor}\leq \frac{deh}{p-1}.
    \]
\end{lem}

\begin{proof}
    The proof of \cite[Lemma 4.1]{Hel12} works verbatim here, except that we should take into account the height $h$, which is equal to $1$ in \cite[Lemma 4.1]{Hel12}. For the reader's convenience, we copy the proof here.
    
Note that there is a filtration $0=\calF^{n+1}\subseteq\calF^n\subseteq\cdots \subseteq \calF^0=\frakM/p\frakM$ such that 
\[
\calF^i/\calF^{i+1}=(\frakM\cap p^iM)/(p\frakM\cap p^{i}M+\frakM\cap p^{i+1}M).
\]
In particular, $\calF^0/\calF^1$ is the free part of the quotient $\frakM/p\frakM$ and each $\calF^i/\calF^{i+1}$ for $i\geq 1$ accounts for the $u^\infty$-torsion contributed by $\frakM\cap (p^iM\backslash p^{i+1}M)$.

Note that in the free $k((u))$-module $p^iM/p^{i+1}M$, there are two $k[[u]]$-lattices: \[(p^i\frakM+(\frakM\cap p^{i+1}M))/(\frakM\cap p^{i+1}M)\subseteq(\frakM\cap p^iM)/(\frakM\cap p^{i+1}M)\]
which is induced from the inclusion $p^i\frakM\subseteq \frakM\cap p^iM$.

Under the natural isomorphism induced by multiplication by $p$: $p^iM/p^{i+1}M\xrightarrow{\cong} p^{i+1}M/p^{i+2}M$, the above two lattices are turned into
\[
(p^{i+1}\frakM+(p\frakM\cap p^{i+2}M))/(p\frakM\cap p^{i+2}M)\subseteq (p\frakM\cap p^{i+1}M)/(p\frakM\cap p^{i+2}M)
\]
which is induced from the inclusion $p^{i+1}\frakM\subseteq p\frakM\cap p^{i+1}M$. Note that multiplication by $p$ induces a surjective map $\frakM\cap p^iM\twoheadrightarrow p\frakM\cap p^{i+1}M$. In fact for any $x\in p\frakM\cap p^{i+1}M$, we have $x=pa=p^{i+1}b$ for some $a\in \frakM$ and $b\in M$. Then $p(a-p^ib)=0\in M$, which means $a-p^ib=p^nc$ for some $c\in M$. As $i\leq n-1$, we see $a\in p^iM\cap \frakM$.

By the theory of elementary divisors, we then have the following equality
\begin{equation*}
    \begin{split}
       & \dim_{k}(\frakM\cap p^iM)/(p^i\frakM+\frakM\cap p^{i+1}M)+\dim_k(\frakM\cap p^{i+1}M)/(p\frakM\cap p^{i+1}M+\frakM\cap p^{i+2}M)\\
       =&\dim_k(\frakM\cap p^{i+1}M)/(p^{i+1}\frakM+\frakM\cap p^{i+2}M)
    \end{split}
\end{equation*}

By using this interpretation, we see
\[\dim_k(\frakM/p\frakM)_{\tor}=\sum_{i=1}^{n}\dim_k\calF^i/\calF^{i+1}=\dim_k(\frakM\cap p^{n}M)/(p^{n}\frakM)
\]
Now we claim that 
\[
u^{\lfloor\frac{he}{p-1}\rfloor}(\frakM\cap p^nM)\subseteq p^n\frakM
\]
In fact, let $j$ be the smallest integer such that
\[
u^{ j}(\frakM\cap p^nM)\subseteq p^n\frakM.
\]
Consider the short exact sequence
\[
0\to p^n\frakM\to \frakM\cap p^nM\to N\to 0.
\]
This means $N=\bigoplus_ik[u]/u^{t_i}$ where $t_i\leq j$ for all $i$ and there exists an $i$ such that $t_i=j$. As the Frobenius $\varphi:W(k)/p^{n+1}[[u]]\to W(k)/p^{n+1}[[u]]$ is flat, we get a short exact sequence
\[
0\to \varphi^*(p^n\frakM)\to \varphi^*(\frakM\cap p^nM)\to \varphi^*(N)\to 0.
\]
Then $\varphi^*(N)=\bigoplus_ik[u]/u^{pt_i}$. So $pj$ is the minimal integer $r$ such that $u^r\varphi^*(\frakM\cap p^nM)\subseteq \varphi^*(p^n\frakM)$.

But we have
\[
p^n\varphi^*(\frakM)\supseteq E^hp^n\frakM=u^{eh}p^n\frakM\supseteq u^{eh+j}(\frakM\cap p^nM)\supseteq u^{eh+j}\varphi^*(\frakM\cap p^nM),
\]
which implies $eh+j\geq pj$. So $j\leq \frac{eh}{p-1}$.

As $\frakM\cap p^nM$ is a lattice in $p^nM$, it is finite free of rank $d$. Hence we get
\[
\dim_{k}(\frakM/p\frakM)_{\tor}\leq \frac{deh}{p-1}.
\]

\end{proof}

We also record the following result.
\begin{lem}[{\cite[Lemma 5.3.14]{EG21}}]\label{finite model}
    Let $(M,\varphi_M)$ be a finite free \'etale $\varphi$-module over $W(k)/p^n((u))$. Then there are only finitely many finitely generated $u$-torsion free Breuil--Kisin modules over $W(k)/p^n[[u]]$ of height $h$, which generate $(M,\varphi_M)$.
\end{lem}

\begin{proof}
    Let $\frakM$ be such a Breuil--Kisin module of height $h$. We may even assume $\frakM$ is finite free up to enlarging $h$. As proved in \cite[Lemma 5.3.14]{EG21},  there exists an integer $i$ such that for any other Breuil--Kisin module $\frakM'$ of height $h$, we have
    \[
    u^{-i}\frakM\subset \frakM'\subset u^i\frakM.
    \]
    Then there are only finitely many such $\frakM'$.
\end{proof}

\begin{proof}[Proof of Proposition \ref{crystalline}]
 Let $\rho_i$ be the induced Galois representation $G_K\to \GL_d(S_{\bar\rho}^{\Box})\to \GL_d(\calO_B/p^i)$ which corresponds to a Laurent $F$-crystal $\bM_i$ over $(\calO_{\Prism}\otimes_{\bZ_p}\calO_B)/p^i[1/\calI_{\Prism}]$. More explicitly, $\bM_i$ corresponds to a stratification $(\calM_i,\epsilon_i)$ where $M_i$ is an \'etale $\varphi$-module over $(W(k)\otimes_{\bZ_p}\calO_B)/p^i((u))$ and $\epsilon_i$ is an isomorphism of $\frakS^{(1)}/p^i[1/E]$-modules
\[
\epsilon_i: \calM_i\otimes_{\frakS/p^i[1/E],p_0}\frakS^{(1)}/p^i[1/E]\cong \calM_i\otimes_{\frakS/p^i[1/E],p_1}\frakS^{(1)}/p^i[1/E].
\]

And $\{\bV_i\}_i$ forms a compatible system in the sense that $\bV_i\otimes_{\bZ/p^i}\bZ/p^{i-1}\cong \bV_{i-1}$. In particular, $(\calM_i,\epsilon_i)\otimes_{\bZ/p^i}\bZ/p^{i-1}\cong (\calM_{i-1},\epsilon_{i-1})$. And by taking the inverse limit, we get a Laurent $F$-crystal $\bM$ over $(\calO_{\Prism}\otimes_{\bZ_p}\calO_B)[1/\calI_{\Prism}]^{\wedge}_p$, which corresponds to a stratification $(\calM,\epsilon)$. This is just the Laurent $F$-crystal attached to the $\calO_B$-representation $\rho$.

Note that each $\rho_i$ only admits finitely many finitely generated prismatic model by Lemma \ref{finite model}. Fix a finitely generated prismatic model $(\frakM_i,\epsilon_i)$ for each $(\calM_i,\epsilon_i)$.  

Then $((\frakM_i/p)/ (\frakM_i/p[u^{\infty}]),\epsilon_i)$ is also a finitely generated prismatic model of $\rho_1$. Let $(\frakM_1',\epsilon_1)$ be one of those appearing infinitely many times. Now among those infinitely many $\frakM_i$ such that $(\frakM_i/p)/ (\frakM_i/p[u^{\infty}])=\frakM_1'$, we can consider $(\frakM_i/p^2)/ (\frakM_i/p^2[u^{\infty}])$, which again are finitely generated prismatic model of $\rho_2$. We choose $(\frakM_2',\epsilon_2)$ to be one of those appearing infinitely many times. In particular, we have $\frakM_1'\cong (\frakM'_{2}/p)/(\frakM'_{2}/p[u^{\infty}])$. Continuing this process, we end up with a system of models $(\frakM_i',\epsilon_i)$ such that $\frakM_i'\cong (\frakM'_{i+1}/p^i)/(\frakM'_{i+1}/p^i[u^{\infty}])$. In particular, the transition map $\frakM_{i+1}'\to \frakM_i'$ is surjective.

Now using Lemma \ref{torsion bounded}, we can further choose a system of models $(\frakM_{n_i}',\epsilon_{n_i})$ such that $\frakM_{n_i}'/p\cong \frakM_{n_{i-1}}'/p$. Let $\frakN_1$ denote $\frakM_{n_i}'/p$.

For each $n_i$, we can consider $p\frakM_{n_i}'$, which is a finitely generated $u$-torsion-free Breuil--Kisin module of height $h$. In particular, $p\frakM_{n_i}'\otimes_{\frakS}\frakS[1/u]\cong p\calM_{n_i}$. There is a short exact sequence
\[
0\to p\frakM_{n_i}'/p^2\frakM_{n_i}'\to \frakM_{n_i}'/p^2\frakM_{n_i}'\to \frakM_{n_i}'/p\frakM_{n_i}'\to 0.
\]

By Lemma \ref{torsion bounded}, $p\frakM_{n_i}'/p^2\frakM_{n_i}'$ has only finitely many choices. So we can find a subsequence $\{m_i\}$ of $\{n_i\}_i$ such that $p\frakM_{m_i}'/p^2\frakM_{m_i}'\cong p\frakM_{m_{i-1}}'/p^2\frakM_{m_{i-1}}'$. In particular, this implies $\frakM_{m_i}'/p^2\frakM_{m_i}'\cong \frakM_{m_{i-1}}'/p^2\frakM_{m_{i-1}}'$. Let $\frakN_2$ denote $\frakM_{m_i}'/p^2\frakM_{m_i}'$. Then we have $\frakN_2/p\cong \frakN_1$ and $\frakN_2[1/u]\cong \calM_2$, which is also compatible with the Frobenius structure.

Continuing this process by considering the short exact sequence
\[
0\to p^s\frakM_{m_i}'/p^{s+1}\frakM_{m_i}'\to \frakM_{m_i}'/p^{s+1}\frakM_{m_i}'\to \frakM_{m_i}'/p^{s}\frakM_{m_i}'\to 0,
\]
we will end up with a compatible system of finitely generated Breuil--Kisin modules $\{\frakN_i\}_i$ such that $\frakN_i/p^{i-1}\cong \frakN_{i-1}$ and $
\frakN_i[1/u]\cong \calM_i$.

Now we consider the inverse limit $\frakN:=\varprojlim_i\frakN_i$. We claim it is a finitely generated module over $\frakS$. In fact, let $(e_1,\cdots, e_t)$ be a minimal set of generators of $\frakN_1$. By Nakayama's lemma,  $(e_1,\cdots, e_t)$ can be lifted to a set of generators of $\frakN_i$.

Then we get a compatible system of short exact sequences given by these generators
\[
0\to K_i\to W(k)/p^{i}[[u]]^{\oplus t}\to \frakN_i\to 0. 
\]
Note that $K_i/p^{i-1}\to W(k)/p^{i-1}[[u]]^{\oplus t}\to \frakN_{i-1}\to 0$ is a right exact sequence. So $K_i/p^{i-1}\to K_{i-1}$ is surjective. So is the natural map $K_i\to K_{i-1}$. Then we get a short exact sequence
\[
0\to \varprojlim_i K_i\to \varprojlim_iW(k)/p^{i}[[u]]^{\oplus t}\to \frakN:=\varprojlim_i\frakN_i\to 0.
\]
So $\frakN$ is finitely generated. 

We naturally have a map $\frakN\to \calM$, which might not be injective. Let $\tilde\frakN$ be its image. We now claim that $\tilde\frakN$ is a finitely generated torsion-free Breuil--Kisin module of height $h$. 

It is clear that $\tilde \frakN$ is finitely generated and torsion-free. As $\varphi:\frakS\to \frakS$ is finite flat, we have $\varphi^*\frakN=\varprojlim_i\varphi^*\frakN_i$. This shows that $\tilde\frakN$ inherits a Frobenius map, which is compatible with $\varphi:\calM\to \calM$. 

We claim that $\frakN/p^i\cong \frakN_i$. Note that there is a compatible system of short exact sequence
\[
0\to T_i\to \frakN/p^i\to \frakN_i\to 0.
\]
Again, $T_i\to T_{i-1}$ is surjective. So we get a short exact sequence
\[
0\to \varprojlim_iT_i\to \frakN=\varprojlim_i\frakN/p^i\to \frakN=\varprojlim_i\frakN_i\to 0.
\]
This shows $\varprojlim_iT_i=0$. As each transition map $T_i\to T_{i-1}$ is surjective, we see that $T_i=0$ for each $i$, i.e. $\frakN/p^i=\frakN_i$. Since $\frakN\otimes_{\frakS}\frakS[1/u]^{\wedge}_p$ is finitely generated over the $p$-complete Noetherian ring $\frakS[1/u]^{\wedge}_p$, we know that $\frakN\otimes_{\frakS}\frakS[1/u]^{\wedge}_p$ is $p$-complete, i.e. $\frakN\otimes_{\frakS}\frakS[1/u]^{\wedge}_p=\varprojlim_i\frakN/p^i[1/u]$. As $\frakN/p^i[1/u]\cong \frakN_i[1/u]\cong \calM_i$, we then have $\frakN\otimes_{\frakS}\frakS[1/u]^{\wedge}_p\cong \calM$ even as an \'etale $\varphi$-modules. This forces $\tilde\frakN\otimes_{\frakS}\frakS[1/u]^{\wedge}_p\cong \calM$ as \'etale $\varphi$-modules. Note that $\tilde\frakN[1/u]$ is finite free by the proof of \cite[Proposition 4.3]{BMS18}. Then the natural map $\varphi^*\tilde \frakN\to \varphi^*\calM$  is injective, which implies $\varphi^*\tilde\frakN\to \tilde\frakN$ is also injective. 

Now we claim the cokernel of the map $\varphi^*\frakN\to \frakN$ is killed by $E^h$. To see this, we note that the cokernel of the map $\varphi^*\frakN_i\to \frakN_i$ is killed by $E^h$. This is because there exists a finitely generated $u$-torsion-free Breuil--Kisin module $\frakM'$ of height $h$ such that $\frakM'/p^i\cong \frakN_i$ by our construction of $\frakN_i$. Then we have

\begin{equation*}
    \xymatrix{
    \varphi^*\frakM'\ar[r]\ar[d] & \frakM'\ar[r]\ar[d] & \bar\frakM'\ar[d]\ar[r]& 0\ar[d]\\
    \varphi^*\frakN_i\ar[r]&\frakN_i\ar[r]&\bar \frakN_i\ar[r]& 0.
    }
\end{equation*}
In particular, the map $\bar\frakM'\to \bar\frakN_i$ is surjective. Since $\bar\frakM'$ is an $\frakS/E^h$-module, we also have $\bar\frakN_i$ is an $\frakS/E^h$-module. We claim that the $\Coker(\varphi^*\frakN\to \frakN)$ is given by $\varprojlim_i\bar\frakN_i$, which is necessarily an $\frakS/E^h$-module. In fact, we have\[
(\varphi^*\frakN\to \frakN)=\varprojlim_i(\varphi^*\frakN_i\to \frakN_i)=\varprojlim_i(\varphi^*\frakN/p^i\to \frakN/p^i).\]

Note that $\Coker(\varphi^*\frakN\to \frakN)/p^i=\Coker(\varphi^*\frakN/p^i\to \frakN/p^i)=\bar\frakN_i$. As $\Coker(\varphi^*\frakN\to \frakN)$ is $p$-complete, we conclude $\Coker(\varphi^*\frakN\to \frakN)=\varprojlim_i\bar\frakN_i$. Again, we have  a surjective map from $\Coker(\varphi^*\frakN\to \frakN)$ to $\Coker(\varphi^*\tilde\frakN\to \tilde\frakN)$. This means $\Coker(\varphi^*\tilde\frakN\to \tilde\frakN)$ is also an $\frakS/E^h$-module. So $\tilde\frakN$ is a finitely generated torsion-free Breuil--Kisin module of height $h$ such that $\tilde \frakN[1/u]^{\wedge}_p\cong \calM$ as \'etale $\varphi$-modules.

In order to get a stratification, we also need to investigate $\varprojlim_i(\frakN_{i}\otimes_{\frakS,p_i}\frakS^{(1)})$. Recall that $\frakN\otimes_{\frakS,p_i}\frakS^{(1)}$ is also $p$-adic complete. So 
\[
\frakN\otimes_{\frakS,p_i}\frakS^{(1)}\cong \varprojlim_j\frakN/p^{j}\otimes_{\frakS,p_i}\frakS^{(1)}\cong \varprojlim_j\frakN_{j}\otimes_{\frakS,p_i}\frakS^{(1)}
\]
for $i=0,1$.

Then we get a stratification $(\frakN,\epsilon: \frakN\otimes_{\frakS,p_0}\frakS^{(1)}\cong\frakN\otimes_{\frakS,p_1}\frakS^{(1)})$, which is compatible with $(\calM,\epsilon)$. This gives us a stratification $(\tilde\frakN,\epsilon: \tilde \frakN\otimes_{\frakS,p_0}\frakS^{(1)}\cong \tilde\frakN\otimes_{\frakS,p_1}\frakS^{(1)})$ by taking the image as $p_i:\frakS\to \frakS^{(1)}$ is (faithfully) flat. By Lemma \ref{observation}, we conclude $\rho$ is a crystalline representation with Hodge--Tate weights in $[0,h]$.

\end{proof}

\begin{thm}\label{main2}
    The $\bZ_p$-flat quotient of $S^{h,\Box}_{\bar \rho}$ is exactly Kisin's crystalline deformation ring, i.e. $S^{h,\Box}_{\bar \rho}/(S^{h,\Box}_{\bar \rho}[p^\infty])=R_{\bar \rho}^{\cris,h,\Box}$.
\end{thm}
\begin{proof}
    Note that Kisin's crystalline deformation ring is uniquely characterised as the $\bZ_p$-flat quotient of $R_{\bar \rho}^{\Box}$ such that for any finite flat $W(\bF)$-algebra $\calO_B$, a map $R_{\bar\rho}^{\Box}\to \calO_B$ factors through $R_{\bar \rho}^{\cris,h,\Box}$ if and only if the induced representation $G_K\to \GL_n(\calO_B)$ is crystalline with Hodge--Tate weights in $[0,h]$. Then this theorem follows from Propositions \ref{crys-factorisation} and \ref{crystalline}.
\end{proof}

\section{Universal prismatic $F$-crystal}
In \cite{kisin2008potentially}, Kisin obtained a finitely generated Breuil--Kisin module with coefficients in the deformation ring of finite height $h$, denoted $A^{\leq h}$ in \textit{loc.cit}. As Kisin's crystalline deformation ring $R_{\bar \rho}^{\cris,h,\Box}$ is a quotient of $A^{\leq h}$, there is an induced finitely generated Breuil--Kisin module $\frakM_{Kis}$ of height $h$ over the crystalline deformation ring, which plays a key role in his investigation of the properties of $R_{\bar \rho}^{\cris,h,\Box}$.

In this section, we aim to put more structures on this Breuil--Kisin module. Indeed, we will equip it with a stratification, which makes it a finitely generated prismatic $F$-crystal. Its \'etale realization corresponds to the universal crystalline Galois deformation $G_K\to \GL_n(R_{\bar \rho}^{\cris,h,\Box})$.

We write $R=R_{\bar \rho}^{\cris,h,\Box}$ with maximal ideal $\frakm$ for simplicity in this section. For each $r>0$, Theorem \ref{main2} shows the corresponding Galois deformation $G_K\to \GL_d(R/\frakm^r)$ admits a finitely generated prismatic $F$-model. 

Let us first set up the notations. Let $\frakS^{(i)}_R:=\varprojlim_r(\frakS^{(i)}\otimes_{\bZ_p}R/\frakm^r)$. In particular, $\frakS_R=(W(k)\otimes_{\bZ_p}R)[[u]]$. Note that the coface maps $\frakS_R\to \frakS_R^{(i)}$ are all flat by \cite[Lemma 3.5]{DLMS24}. Attached to the Galois representation $G_K\to \GL_d(R)$, we have the finite free Laurent $F$-crystal $(N,\varphi,\epsilon_N: N\otimes_{\frakS_R,p_0}\frakS^{(1)}_R\cong N\otimes_{\frakS_R,p_1}\frakS^{(1)}_R$, where $(N,\varphi)$ is a finite free \'etale $\varphi$-module over $\frakS_R[1/u]^\wedge_\frakm$ of rank $n$. We simply call $N$ the Laurent $F$-crystal and write $\Phi: \varphi^*N\xrightarrow{\cong} N$.

\begin{thm}\label{universal prismatic crystal}
 There is a finitely generated prismatic $F$-model $(\frakM,\epsilon:\frakM\widehat\otimes_{\frakS_R,p_0}\frakS^{(1)}_R\cong \frakM\widehat\otimes_{\frakS_R,p_1}\frakS^{(1)}_R)$  of height $h$ , where the completion is $\frakm$-adic completion, whose \'etale realization gives the universal deformation $G_K\to \GL_d(R)$.
\end{thm}

The method of proving Proposition \ref{crystalline} does not work here as there is no uniform bound as in Lemma \ref{torsion bounded} for $R/\frakm^r$. We will use a different method to prove Theorem \ref{universal prismatic crystal}.

The idea is to first prove the existence of a prismatic $F$-model for $N/p^r$ for each $r$. As $R/p^r$ is not a finite Artinian ring, we do not have the existence of a model a priori. We can try to take the inverse limit along $n$ of models over $R/{(p^r,\frakm^n)}$. But this essentially requires some boundedness results. The existence of the Breuil--Kisin module in \cite{kisin2008potentially} is the key here. 

\begin{lem}\label{R-lattice}
    There are isomorphisms of \'etale $\varphi$-modules over $\frakS_R[1/u]^\wedge_{\frakm}[1/p]$
    \[
    \frakM_{Kis}\otimes_{\frakS_R}\frakS_R[1/u]^\wedge_{\frakm}[1/p]\cong N[1/p].
    \]
\end{lem}
\begin{proof}
We use some notations from \cite[Section 1.7]{kisin2008potentially}. In particular, we let $A$ be the framed unrestricted deformation ring and $A^{\leq h}$ be the deformation ring of finite height $h$. Note that the crystalline deformation ring $R$ is a quotient of $A^{\leq h}$.  The map $\Theta_{\frakS_{A^{\leq h}}}:X_{A^{\leq h}}:=\Spec(\calL_A^h)\times_{\Spec({A^{\leq h}})}\Spec(\frakS_{A^{\leq h}})$ is proper and is an isomorphism after inverting $p$. There is a universal sheaf $\underline\frakM$ on $X_{A^{\leq h}}$ and $\frakM_{Kis}=H^0(X_{A^{\leq h}}\times_{\frakS_{A^{\leq h}}}\frakS_{R},\underline\frakM)$. By the definition of $\Spec(\calL_A^{\leq h})$, we see that $\underline\frakM\otimes_{\frakS_{A^{\leq h}}}\frakS_R/\frakm^n[1/u]=\Theta^*N$ for each $n$, where $\Theta$ is the pullback of $\Theta_{\frakS_{A^{\leq h}}}$ along the composite $\Spec(\frakS_R[1/u]^\wedge_\frakm)\to\Spec(\frakS_R)\to \Spec(\frakS_{A^{\leq h}})$.

By \cite[Proposition 1.6.4]{kisin2008potentially}, we see $\Theta_{\frakS_{A^{\leq h}}}[1/p]$ is an isomorphism. So
\[
\frakM_{Kis}\otimes_{\frakS_R}\frakS_R[1/u]^\wedge_{\frakm}[1/p] \cong H^0(X_{A^{\leq h}}\times_{\frakS_{A^{\leq h}}}\frakS_R[1/u]^\wedge_{\frakm}[1/p],\Theta^*N)\cong N[1/p].
\]

\end{proof}

The significance of the existence of $\frakM_{Kis}$ is that it enables us to prove  the following analogue of \cite[Lemma 5.3.14]{EG21} for coefficients rings $R/p^r$.
\begin{lem}\label{bound BK-lattice}
    Fix an integer $r>0$. There exists a finitely generated Breuil--Kisin module $\Lambda_r\subseteq N_r=N/p^r$ of height $h$ over $\frakS_R/p^r$ and an integer $c_r\geq 0$, which only depends on $r$, such that for any Breuil--Kisin lattice $\frakM_{r,n}$ of height $h$ in $N_{r,n}:=N/{(p^r,\frakm^n)}$, we have 
    \[
    u^{c_r}\Lambda_{r,n}\subseteq \frakM_{r,n}\subseteq u^{-c_r}\Lambda_{r,n}
    \]
    where $\Lambda_{r,n}=\Ima(\Lambda_r\to N_{r,n})$. In particular, for any Breuil--Kisin lattice $\frakM_r$ in $N_r$, we also have
    \[
    u^{c_r}\Lambda_{r}\subseteq \frakM_{r}\subseteq u^{-c_r}\Lambda_{r}.
    \]
\end{lem}

\begin{proof}
By Lemma \ref{R-lattice}, there exist a finitely generated Breuil--Kisin module $\tilde\frakM$ of height $h$ over $\frakS_R$ and an integer $D$ such that 
\[
p^DN\subseteq \tilde\frakM\otimes\frakS_{R}[1/u]^\wedge_\frakm\subseteq N
\]

    By the proof of \cite[Lemma 5.3.14]{EG21}, we know that for any Breuil--Kisin modules $L_1,L_2$ of height $h$ satisifying $L_1[1/u]\subseteq L_2[1/u]\subseteq N_{s,n}$, there exists an integer $t_s$ which only depends on $s, h$ and $E(u)$ such that
    \[
    L_1\subseteq u^{-t_s}L_2.
    \]

Now we let $L_{s,n}$ be the maximal finitely generated Breuil--Kisin lattice of height $h$ in $N_{s,n}$. Put $L_s:=\varprojlim_nL_{s,n}\subseteq N_s$. Consider $\tilde\frakM_{s,n}:=\Ima(\tilde\frakM\to N_{s,n})$. Then we have
\[
p^DL_{s,n}\subseteq u^{-t_s}\tilde\frakM_{s,n}.
\]
Taking the inverse limit along $n$, we get
\[
p^DL_s\subseteq u^{-t_s}\tilde\frakM_s
\]
where $\tilde\frakM_s:=\Ima(\tilde\frakM\to N_s)$. Note that each $\tilde\frakM_{s,n}$ is finite. By Lemma \ref{nonempty}, we have $\tilde\frakM_s=\varprojlim_n\tilde\frakM_{s,n}$.

Now we get each $p^DL_s$ is a finitely generated over $\frakS_R/p^s$ of height $h$. Put $\Lambda_r:=\Ima(L_{r+D}\to N_{r})$. To see $\Lambda_r$ is finitely generated over $\frakS_R/p^s$, we consider the injection $\iota_D:N/p^r\xrightarrow{\times p^D}N/{p^{r+D}N}$. Note that $\iota_D(\Lambda_r)=p^DL_{r+D}$. So $\Lambda_r$ is finitely generated. In particular, $\Lambda_r$ is a finitely generated Breuil--Kisin module of height $h$ over $\frakS_R/p^r$.

It remains to show the uniform bound. Note that for each $m\geq n$, the $\Ima(L_{s,m}\to N_{s,n})$ is a model of height $h$ in $N_{s,n}$. So we get
\[
u^{t_s}L_{s,n}\subseteq \Ima(L_{s,m}\to N_{s,n}).
\]
Taking the inverse limit along $m$, we get
\[
u^{t_s}L_{s,n}\subseteq \Ima(L_{s}\to N_{s,n}).
\]

So $\Lambda_{r,n}$ is a Breuil--Kisin lattice of height $h$ in $N_{r,n}$. By \cite[Lemma 5.3.14]{EG21}, there exists an integer $c_r$ such that for any Breuil--Kisin lattice $\frakM_{r,n}$ of height $h$ in $N_{r,n}:=N/{(p^r,\frakm^n)}$, we have 
    \[
    u^{c_r}\Lambda_{r,n}\subseteq \frakM_{r,n}\subseteq u^{-c_r}\Lambda_{r,n}.
    \]

   For any Breuil--Kisin lattice $\frakM_r$ of height $h$ in $N_r$, we consider its image $\frakM_{r,n}:=\Ima(\frakM_r\to N_{r,n})$ for each $n$. By Lemma \ref{nonempty} $\frakM_r=\varprojlim_n\frakM_{r,n}$. So we also have 
   \[
       u^{c_r}\Lambda_{r}\subseteq \frakM_{r}\subseteq u^{-c_r}\Lambda_{r}.
   \]
\end{proof}

\begin{lem}\label{nonempty}
    Let $A$ be a Noetherian ring which is $I$-complete for some ideal $I$ and $A/I^n$ is finite for all $n$. Let $M$ be a finitely generated $A$-module equipped with $I$-adic topology. Suppose there is an inverse system of surjective continuous maps $\{M\twoheadrightarrow M_n\}_n$ where each $M_n$ is finite with discrete topology. Then the induced map $M\to \varprojlim_nM_n$ is surjective.
\end{lem}
\begin{proof}
    It is clear that $M$ is a compact Hausdorff space. Let $\{x_n\}_n\in \varprojlim_nM_n$. Put $M(x_n)$ to be the set of preimages of $x_n$ under the surjective map $M\to M_n$. Then each $M_n$ is closed in $M$ and $M_{n+1}\subseteq M_n$. So $\bigcap_nM(x_n)\neq \emptyset$. We are done.
\end{proof}

So far, we only made preparations about the Breuil--Kisin modules. In order to get a prismatic $F$-model, we also need to construct stratifications. The following theorem due to Chevalley will be crucial in this aspect.

\begin{thm}[Chevalley's theorem]\label{Chevalley}
    Let $(A,\frakm)$ be a complete Noetherian local ring. Let $M$ be a finite $A$-module and $H_1\supseteq H_2\cdots$ be submodules of $M$. Write $H=\cap_i H_i$. Then for every $a\geq 1$, there is an integer $i(a)$ such that
    \[
    H_{i(a)}\subseteq H+\frakm^aM.
    \]
\end{thm}
\begin{proof}
Just apply \cite[Theorem 19.9]{Hochster} to $M/H$.
\end{proof}

\begin{lem}\label{tenstor product commutes with cap}
    Let $p_j:\frakS_R\to \frakS_R^{(i)}$ be any coface map. Suppose $T_1\supseteq T_2\cdots$ are $\frakS_R$-submodules of $N_r$ such that 
    \[
    u^c\Lambda_r\subseteq T_i\subseteq u^{-c}\Lambda_r
    \]
    for some constant $c>0$. Then we have
    \[
\frakS_R^{(i)}\otimes_{p_j,\frakS_R}(\bigcap_iT_i)=\bigcap_i(\frakS_R^{(i)}\otimes_{p_j,\frakS_R}T_i).
    \]
\end{lem}
\begin{proof}
    Write $M_1=u^{-c}\Lambda_r$, $M_2=u^c\Lambda_r$ and $T=\bigcap_i T_i$. Then $M_1/M_2$ is finite over the complete Noetherian local ring $R/p^r$. Applying Theorem \ref{Chevalley} to $T_i/M_2\subseteq M_1/M_2$ shows for each $a$, we have
    \[
    T_{i(a)}\subseteq T+\frakm^aM_1.
    \]
    As $p_i:\frakS_R\to \frakS_R^{(i)}$ is flat, we have
    \[
 \bigcap_i(\frakS_R^{(i)}\otimes_{p_j,\frakS_R}T_i)\subseteq \frakS_R^{(i)}\otimes_{p_i,\frakS_R}T+\frakm^a(\frakS_R^{(i)}\otimes_{p_i,\frakS_R}M_1)
    \]
    for all $a$. Hence the image of $\bigcap_i(\frakS_R^{(i)}\otimes_{p_j,\frakS_R}T_i)$ under the map $\frakS_R^{(i)}\otimes M_1\to \frakS_R^{(i)}\otimes M_1/T$ is $0$ as $\frakS_R^{(i)}\otimes M_1$ is $\frakm$-adic separated. This gives
    \[
    \bigcap_i(\frakS_R^{(i)}\otimes_{p_j,\frakS_R}T_i)\subseteq \frakS_R^{(i)}\otimes_{p_j,\frakS_R}(\bigcap_iT_i).
    \]
    The reverse inclusion is obvious. So we are done.
\end{proof}

Now we are ready to prove Theorem \ref{universal prismatic crystal}.

\begin{proof}[Proof of Theorem \ref{universal prismatic crystal}]
    Under our assumption, each $N_{r,n}$ admits a finitely generated prismatic $F$-model of height $h$. Let $\frakM_{r,n}$ be the sum of  all such models in $N_{r,n}$. Here $\frakM_{r,n}$ is short for the stratification $(\frakM_{r,n},\epsilon_{r,n}: \frakM_{r,n}\otimes_{\frakS_R,p_0}\frakS_R^{(1)}\cong \frakM_{r,n}\otimes_{\frakS_R,p_1}\frakS_R^{(i)}$. By Lemma \ref{bound BK-lattice}, the sum of all underlying Breuil--Kisin lattices of height $h$ is contained in $u^{-t_r}\Lambda_{r,n}$. So $\frakM_{r,n}$ is a finite sum of finitely generated prismatic $F$-models, which itself is a finitely generated prismatic $F$-model of height $h$. This means $\frakM_{r,n}$ is actually the maximal finitely generated prismatic $F$-model of height $h$ in $N_{r,n}$. The reduction map $N_{r,n+1}\to N_{r,n}$ restricts to a map $\frakM_{r,n+1}\to \frakM_{r,n}$.

    Now we consider $\frakM_r:=\varprojlim_n\frakM_{r,n}$. Lemma \ref{bound BK-lattice} gives $u^{t_r}\Lambda_{r,n}\subseteq \frakM_{r,n}\subseteq u^{-t_r}\Lambda_{r,n}$. Taking inverse limit with respect to $n$, we see $u^{t_r}\Lambda_{r}\subseteq \frakM_{r}\subseteq u^{-t_r}\Lambda_{r}$. As $\varphi:\frakS_R\to \frakS_R$ is finite free, we have $\varphi^*\frakM_r=\varprojlim_n\varphi^*\frakM_{r,n}$.  Note that $E$ is a non-zero-divisor in each $\frakS_R/(p^r,\frakm^n)=(R/(p^r,\frakm^n)\otimes_{\bZ/p^r}W_r(k))[[u]]$. So $\varprojlim_nE^h\frakM_{r,n}=E^h\frakM_r$. This implies $\frakM_r$ is a finitely generated Breuil--Kisin module of height $h$ in $N_r$ and $\frakM_r\otimes_{\frakS_R}\frakS_R[1/u]^\wedge_\frakm\cong N_r$.

    We claim $\frakM_r$ together with the restriction of $\epsilon_{N_r}$ forms a stratification.

The finitely generated $\frakS_R$-module $u^{-t_r}\Lambda_r$ contains $\frakM_r$ and contains all $\frakM_{r,n}$ after reduction. Consider the stratification $\epsilon_{N_r}:N_r\otimes_{\frakS_R,p_0}\frakS_R^{(1)}\cong N_r\otimes_{\frakS_R,p_1}\frakS_R^{(1)}$. We claim that the image $\epsilon_{N_r}(u^{-t_r}\Lambda_r\otimes_{\frakS_R,p_0}\frakS_R^{(1)})$ is contained in $u^{-t_r-\alpha_r}\Lambda_r\otimes_{\frakS_R,p_1}\frakS_R^{(1)}$ for some integer $\alpha_r$.

Note that \[\epsilon_{N_{r,n}}(u^{-t_r}\Lambda_{r,n}\otimes_{p_0}\frakS_R^{(1)})\subseteq p_0(u)^{-2t_r}\epsilon_{N_{r,n}}(\frakM_{r,n}\otimes_{p_0}\frakS_R^{(1)})=p_0(u)^{-2t_r}\frakM_{r,n}\otimes_{p_1}\frakS_R^{(1)}\] and 
\[
p_0(u)^{-2t_r}\frakM_{r,n}\otimes_{p_1}\frakS_R^{(1)}\subseteq p_0(u)^{-2t_r}p_1(u)^{-t_r}\Lambda_{r,n}\otimes_{p_1}\frakS_R^{(1)}
\]
 Since $(p_0(E))=(p_1(E))\subseteq \frakS^{(1)}$ and $E\equiv u^e\mod p$, we have $p_1(u)^e\in (p_0(u)^e,p)\frakS^{(1)}/p^r$ and hence $p_1(u)^{er}\in p_0(u)^e\frakS^{(1)}/p^r$. Put
 \[
 \alpha_r=er\lceil \frac{2t_r}{e}\rceil
 \]
so that $p_1(u)^{\alpha_r}\in p_0(u)^{2t_r}\frakS^{(1)}/p^r$ and 
$\epsilon_{N_{r,n}}(u^{-t_r}\Lambda_{r,n}\otimes_{p_0}\frakS_R^{(1)})\subseteq p_1(u)^{-t_r-\alpha_r}\Lambda_{r,n}\otimes_{p_1}\frakS_R^{(1)}$.

Now we claim that $\Lambda_r\otimes_{\frakS_R,p_i}\frakS_R^{(1)}\cong \varprojlim_n(\Lambda_{r,n}\otimes_{\frakS_R,p_i}\frakS_R^{(1)})$.
As there are inclusions $\Lambda_r\otimes_{\frakS_R,p_i}\frakS_R^{(1)}\subseteq N_r\otimes_{\frakS_R,p_i}\frakS_R^{(1)}$ and $\varprojlim_n(\Lambda_{r,n}\otimes_{\frakS_R,p_i}\frakS_R^{(1)})\subseteq N_r\otimes_{\frakS_R,p_i}\frakS^{(1)}$, we then have an inclusion $\Lambda_r\otimes_{\frakS_R,p_i}\frakS_R^{(1)}\subseteq \varprojlim_n(\Lambda_{r,n}\otimes_{\frakS_R,p_i}\frakS_R^{(1)})$. To see the surjectivity, let $K_n=\Ker(\Lambda_r\to \Lambda_{r,n}$. Then $\bigcap_n K_n=0$. By Theorem \ref{Chevalley}, we see for each $a$, there exists an integer $n(a)$ such that $K_{n(a)}\subseteq (\frakm,u)^a\Lambda_r$. As $p_i:\frakS_R\to \frakS_R^{(1)}$ is flat, we see $\Lambda_r\otimes_{\frakS_R,p_i}\frakS_R^{(1)}$ is $(\frakm,u)$-adic complete and $\Ker(\Lambda_r\otimes_{\frakS_R,p_i}\frakS_R^{(1)}\to \Lambda_{r,n}\otimes_{\frakS_R,p_i}\frakS_R^{(1)})=K_n\otimes_{\frakS_R,p_i}\frakS_R^{(1)}$. For any $\{x_n\}_n\in \varprojlim_n\Lambda_{r,n}\otimes_{\frakS_R,p_i}\frakS_R^{(1)}$, we can choose a lift $\tilde x_n$ for each $x_n$ and then $\{\tilde x_n\}_n$ is a Cauchy sequence in $\Lambda_{r}\otimes_{\frakS_R,p_i}\frakS_R^{(1)}$. So we get the desired surjectivity.

Write $\tilde\Lambda_r=u^{-t_r-\alpha_r}\Lambda_r$ and let $T_n\subseteq \tilde\Lambda_r$ be the preimage of $\frakM_{r,n}$ under the reduction map $\tilde\Lambda_r\to N_{r,n}$. Then $\bigcap_n T_n=\frakM_r$. For any $x\in \frakM_r\otimes_{\frakS_R,p_0} \frakS_R^{(1)}$, the reduction of $\epsilon_{N_r}(x)$ lies in $\frakM_{r,n}\otimes_{\frakS_R,p_1}\frakS_R^{(1)}$ for each $n$. Again as $p_i:\frakS_R\to \frakS_R^{(1)}$ is flat, we have $\epsilon_{N_r}(x)\in T_n\otimes_{\frakS_R,p_1}\frakS_R^{(1)}$ for each $n$. Lemma \ref{tenstor product commutes with cap} then shows $\epsilon_{N_r}(x)\in \frakM_r\otimes_{\frakS_R,p_1}\frakS_R^{(1)}$. Apply the same argument to $\epsilon^{-1}$ and then we get the isomorphism $\frakM_r\otimes_{\frakS_R,p_0}\frakS_R^{(1)}\cong \frakM_r\otimes_{\frakS_R,p_1}\frakS_R^{(1)}$. The prismatic $F$-model $(\frakM_r,\epsilon_r)$ is maximal in the sense that any other finitely generated prismatic $F$-model of height $h$ must be contained in $\frakM_r$ as their reductions must be contained in $\frakM_{r,n}$ for each $n$.

Now consider $\frakM:=\varprojlim_r\frakM_r\subseteq N$. We claim that $\frakM$ is a finitely generated Breuil--Kisin lattice over $\frakS_R$ of height $h$ in $N$.

Note that multiplication by $p$ gives an isomorphism $\iota_r:N_r\xrightarrow{\times p}pN_{r+1}$. The maximality of $\frakM_r$ then implies 
\[
\frakM_{r+1}\cap pN_{r+1}=\iota_r(\frakM_r).
\]
Indeed, $\frakM_{r+1}\cap pN_{r+1}$ is a prismatic $F$-model of $pN_{r+1}$, hence must be contained in $\iota_r(\frakM_r)$. Conversely,  $\iota_r(\frakM_r)$ is also a prismatic $F$-model of $pN_{r+1}$, hence is contained in $\frakM_{r+1}$.

It then follows 
\begin{equation}\label{saturation}
\frakM\cap p^rN=p^r\frakM
\end{equation}
for every $r$. Consider the continuous map $N\to N/p^r$ where $N$ is equipped with $p$-adic topology and $M/p^r$ is discrete. Then $\frakM$ is the intersection of the preimages of the closed subsets $\frakM_{r}\subseteq N_r$, hence is closed in $N$. Then we must have $\frakM=\varprojlim_r\frakM/p^r$. 
Equation \ref{saturation} implies the injection $\frakM/p\hookrightarrow \frakM_1\subseteq N_1$. So $\frakM/p$ is finite over $\frakS_R/p$. Nakayama's lemma then shows $\frakM$ is finite over $\frakS_R$. Again as $\varphi:\frakS_R\to \frakS_R$, we see that $\frakM$ is a finitely generated Breuil--Kisin module over $\frakS_R$ of height $h$.

It remains to show the natural map $\tau:\frakM\otimes_{\frakS_R}\frakS_R[1/u]^\wedge_\frakm\to N$ is an isomorphism. To show this, we first prove $\frakM/p\otimes_{\frakS_R}\frakS_R[1/u]^\wedge_\frakm\cong N_1$.

Let $I_{s,r}$ be the image of $\frakM_s$ in $N_r$ for any $s\geq r$. Then each $I_{s,r}$ is a prismatic $F$-model of height $h$ in $N_r$. Therefore we have for every $r$
    \[
    u^{t_r}\Lambda_r\subseteq I_{s,r}.   
    \]
    For any $v\in \bigcap_{s\geq r}I_{s,r}$, we define
    \[
    X_{s,r}(v):=\{x\in \frakM_s: x\mod p^r=v\}.
    \]
    Equip each $\frakM_r$ with the $(\frakm,u)$-adic topology. Note that $\frakS_R/p^r$ is complete with respect to the $(\frakm,u)$-adic topology. Every finite $\frakS_R/p^r$-module is then compact Hausdorff. The reduction maps $\frakM_{s}\to \frakM_{r}$ are continuous, so each $X_{s,r}(v)$ is closed and compact. For a fixed $r$, the inverse limit $\varprojlim_{s\geq r} X_{s,r}(v)$ is then non-empty. So $v$ lifts to $\frakM$. This means $u^{t_r}\Lambda_r\subseteq \bigcap_{s\geq r}I_{s,r}=\frakM/p^r\subseteq N_r$. Then we get
    \[
    \frakM/p^r\otimes_{\frakS_R}\frakS_R[1/u]^\wedge_\frakm\cong N_r.
    \]
    Hence the natural map $\tau:\frakM\otimes_{\frakS_R}\frakS_R[1/u]^\wedge_\frakm\to N$ is an isomorphism modulo $p$. The $\Coker(\tau)$ is a finite $\frakS_R[1/u]^\wedge_\frakm$-module and then vanishes by Nakayama's lemma. As $N$ is finite free, we see $\frakM\otimes_{\frakS_R}\frakS_R[1/u]^\wedge_\frakm\cong N\oplus \Ker(\tau)$. Hence $\Ker(\tau)$ also vanishes by Nakayama's lemma. 

    Now by applying Lemma \ref{tenstor product commutes with cap} to $\{I_{s,r}\}_{s\geq r}$, we then get the stratification
    \[
    \epsilon: \frakM\otimes_{\frakS_R,p_0}\frakS_R^{(1)}\cong \frakM\otimes_{\frakS_R,p_1}\frakS_R^{(1)},
    \]
    which is the restriction of $\epsilon_{N}$.

    For the maximality, note that the mod $p^r$ reduction of any prismatic $F$-model in $N$ must be contained in $\frakM_r$ due to the maximality of the latter, and hence itself must be contained in $\frakM$.
\end{proof}

Next we want to compare the underlying Breuil--Kisin module of the universal prismatic $F$-crystal over Kisin's crystalline deformation ring with the finitely generated Breuil--Kisin module over the deformation ring in \cite{kisin2008potentially}. 

\begin{lem}\label{vanishing}
Let  $M$ be a finite $\frakS_R[1/p]$-module equipped with a semilinear Frobenius map $\varphi:M\to M$ such that $(\varphi^*M)[1/E]\cong M[1/E]$. If $M\otimes_{\frakS_R}\frakS_R[1/u]^\wedge_\frakm=0$, then $M=0$.
\end{lem}
\begin{proof}
    Choose a finite $\frakS_R$-submodule $\frakM\subseteq M$ such that $\frakM[1/p]=M$. As $\frakM$ is $p$-torsion free and $\frakS_R[1/u]^\wedge_\frakm$ is flat over $\frakS_R$, the module $\frakM\otimes_{\frakS_R}\frakS_R[1/u]^\wedge_\frakm$ is also $p$-torsion free. As $\frakM\otimes_{\frakS_R}\frakS_R[1/u]^\wedge_\frakm[1/p]=0$, we see $\frakM\otimes_{\frakS_R}\frakS_R[1/u]^\wedge_\frakm=0$. So $(\frakM/\frakm\frakM)[1/u]=0$.

    This means some power of $u$ kills $\frakM/\frakm\frakM$. Hence $\frakM/\frakm$ is finite over $R/\frakm\otimes_{\bF_p}k$, and hence finite over $R/\frakm$. Nakayama's lemma shows $\frakM$ is finite over $R$. Thus $M$ is finite over $R[1/p]$.

    Define 
    \[
    D(u)=\prod_{\tau\in \Gal(K_0/\bQ_p)}\tau(E(u))\in \bZ_p[u].
    \]
   Note that the polynomial $E$ divides $D$ and every root of $D$ has $p$-adic valuation $1/e$ under the normalization $v(p)=1$.

   Let $F$ be the residue field of any maximal ideal of $R[1/p]$ and put
   \[
   V=M\otimes_{R[1/p]}F.
   \]
   Let $m=\dim_FV$ and $f(X)=\det(X-u\mid V)$.

   As $F$-vector spaces, we have $\varphi^*V:=(\frakS_R\otimes_{\varphi,\frakS_R}M)\otimes_{R[1/p]}F\cong V^{\oplus p}$. Thus $f(X^p)=\det(X-u\mid \varphi^*V)$.

   By the assumption, we have $(\frakS_R\otimes_{\varphi,\frakS_R}M)[1/D]\cong M[1/D]$. So $g(X)=\frac{f(X^p)}{f(X)}$ has zeros and poles only among the roots of $D$.

   Since $D(0)\neq 0$, comparison of multiplicities at $0$ gives $p\cdot ord_0(f)=ord_0(f)$. Therefore $f(0)\neq 0$ and $g(0)=1$.

   Write $g(X)=\prod_{\alpha}(X-\alpha)^{c_\alpha}$ where each $\alpha$ is a root of $D$. Comparing degrees, we get
   \[
   \sum_{\alpha}c_{\alpha}=(p-1)m.
   \]
   It follows that $0=v(g(0))=\sum_{\alpha}c_{\alpha}v(\alpha)=\frac{(p-1)m}{e}$.

   Hence $m=0$. This means the finite $R[1/p]$-module $M$ vanishes at any closed point of $R[1/p]$. So $M=0$.
\end{proof}

\begin{lem}
    Assume $R[1/p]$ is regular. Then for any Breuil--Kisin lattice $\frakM\subseteq N$, the module $\frakM[1/p]$ is finite projective over $\frakS_R[1/p]$.
\end{lem}
\begin{proof}
    The ring $\frakS_R[1/p]$ is regular of finite dimension under the assumption. For any $i>0$, let $Q_i=\Ext^j_{\frakS_R[1/p]}(\frakM[1/p], \frakS_R[1/p])$. Then we have
    \[
    Q_i\otimes_{\frakS_R}\frakS_R[1/u]^\wedge_\frakm\cong \Ext^i_{\frakS_R[1/u]^\wedge_\frakm}(\frakM[1/p]\otimes_{\frakS_R}\frakS_R[1/u]^\wedge_\frakm=N[1/p],\frakS_R[1/u]^\wedge_\frakm[1/p])=0.
    \]
    As the Frobenius map is flat, we get
    \[
    \Ext^i_{\frakS_R[1/p]}(\varphi^*\frakM[1/p],\frakS_R[1/p])=\varphi^*Q_i.
    \]
    So $(\varphi^*Q_i)[1/E]\cong Q_i[1/E]$. By Lemma \ref{vanishing}, we see $Q_i=0$.

    Since $\frakS_R[1/p]$ is regular of finite dimension, $\frakM[1/p]$ has finite projective dimension. The vanishing of all Ext groups then implies $\frakM[1/p]$ is projective.

\end{proof}

\begin{rmk}
   The existence of a prismatic $F$-model $\bM$ of the universal crystalline Galois deformation $G_K\to \GL_d(R_{\bar\rho}^{\cris,\Box})$ has some interesting implications. For example, one can consider the attached Hodge--Tate crystal $\bM/\calI_{\Prism}$ (note that $\bM/\calI_{\Prism}$ is indeed $p$-torsion free). Then by the classification of Hodge--Tate crystals in \cite{gao2023hodge}, one gets a single integral Sen operator $\Theta: \bM/\calI_{\Prism}(\frakS,E))\to \bM/\calI_{\Prism}(\frakS,E))$ which controls the Hodge--Tate weights of the parametrised crystalline representations.
\end{rmk}

\addcontentsline{toc}{section}{References}
\bibliographystyle{alpha}
\bibliography{ref}
\end{document}